\documentclass{article}%
\usepackage{amssymb}
\usepackage{amsfonts}
\usepackage{graphicx}
\usepackage{amsmath}%
\newtheorem{theorem}{Theorem}

\newtheorem{corollary}[theorem]{Corollary}

\newtheorem{lemma}[theorem]{Lemma}

\begin{document}

\title{Bernoulli actions are weakly contained in any free action}
\author{Mikl\'{o}s Ab\'{e}rt and Benjamin Weiss}
\maketitle

\begin{abstract}
Let $\Gamma$ be a countable group and let $f$ be a free probability measure
preserving action of $\Gamma$. We show that all Bernoulli actions of $\Gamma$
are weakly contained in $f$.

It follows that for a finitely generated group $\Gamma$, the cost is maximal
on Bernoulli actions for $\Gamma$ and that all free factors of i.i.d. of
$\Gamma$ have the same cost.

We also show that if $f$ is ergodic, but not strongly ergodic, then $f$ is
weakly equivalent to $f\times I$ where $I$ denotes the trivial action of
$\Gamma$ on the unit interval. This leads to a relative version of the
Glasner-Weiss dichotomy.

\end{abstract}

\section{Introduction}

Let $\Gamma$ be a countable group. By a probability measure preserving
(p.m.p.) action $f$ of $\Gamma$ we mean a representation of $\gamma\in\Gamma$
by $f_{\gamma}$ in the group of measure preserving transformations of
$(X,\mathcal{B},\mu)$. Here $(X,\mathcal{B},\mu)$ is a standard Borel
probability space. The action $f$ is \emph{free}, if for $\mu$-almost every
$x\in X$ and $\gamma\neq\gamma^{\prime}\in\Gamma$ we have $f_{\gamma}(x)\neq
f_{\gamma^{\prime}}(x)$.

Let $f$ and $g$ be p.m.p. actions of $\Gamma$ on $(X,\mathcal{B},\mu)$ and
$(Y,\mathcal{C},\nu)$, respectively. Following the definition of Kechris
\cite{kechbook}, we say that $f$ \emph{weakly contains} $g$ ($f\succeq g$) if
for all Borel subsets $Y_{1},\ldots,Y_{n}\in\mathcal{C}$, finite sets
$S\subseteq\Gamma$ and $\varepsilon>0$ there exist Borel subsets $X_{1}%
,\ldots,X_{n}\in\mathcal{B}$ such that
\[
\left\vert \mu(f_{\gamma}X_{i}\cap X_{j})-\nu(g_{\gamma}Y_{i}\cap
Y_{j})\right\vert <\varepsilon\text{ \ (}1\leq i,j\leq n,\gamma\in S\text{).}%
\]
In particular, the way $g$ acts on finite partitions of the underlying space
can be simulated by $f$ with arbitrarily small error. A special case of weak
containment is when $g$ is a \emph{factor} of $f$, that is, there exists a
surjective measure preserving $\Gamma$-equivariant map from $(X,\mathcal{B}%
,\mu)$ to $(Y,\mathcal{C},\nu)$. We call $f$ and $g$ \emph{weakly equivalent}
if $f\succeq g$ and $g\succeq f$.

Let $\kappa$ be a probability space. The Bernoulli action $\kappa^{\Gamma} $
is defined as the set of maps from $\Gamma$ to $\kappa$, endowed with the
product measure and the shift action by $\Gamma$. These actions are also
called i.i.d. processes on $\Gamma$.

The main result of this paper is the following:

\begin{theorem}
\label{fotetel}Let $\Gamma$ be a countable infinite group and let $f$ be a
free p.m.p. action of $\Gamma$. Then $f$ weakly contains every Bernoulli
action of $\Gamma$. In particular, all free factors of i.i.d.-s of $\Gamma$
are weakly equivalent.
\end{theorem}

Note that one can obtain the original Rokhlin lemma for $\mathbb{Z}$ as a
quick corollary of Theorem \ref{fotetel}. In this language, the lemma says
that every free action of $\mathbb{Z}$ weakly contains the cyclic action of
$\mathbb{Z}$ on $n$ points for all $n>0$. Since weak containment is
transitive, by Theorem \ref{fotetel}, it is enough to prove this for
$\{0,1\}^{\mathbb{Z}}$ with $p(0)=p(1)=1/2$. For a natural number $k$, $0\leq
l\leq n-1$ and $\omega\in\{0,1\}^{\mathbb{Z}}$ let
\[
\left\Vert \omega\right\Vert _{k,l}=%
{\displaystyle\sum\limits_{i=0}^{k-1}}
\omega(in+l)
\]
and let
\[
A_{k}=\left\{  \omega\in\{0,1\}^{\mathbb{Z}}\mid\left\Vert \omega\right\Vert
_{k,0}>\left\Vert \omega\right\Vert _{k,l}+n\text{ for all }1\leq l\leq
n-1\right\}  \text{.}%
\]
Let $T$ denote the shift operator on $\{0,1\}^{\mathbb{Z}}$. Then the sets
$A_{k},TA_{k},\ldots,T^{n-1}A_{k}$ are disjoint, while the measure of $A_{k}$
converges to $1/n$ as $k$ tends to infinity, since for a fixed $k$, the
$\left\Vert \omega\right\Vert _{k,l}$ are independent binomial variables.

As a corollary of Theorem \ref{fotetel}, one gets a new result on cost. The
cost is a numeric invariant of a p.m.p. action of a countable group. Its study
was initiated by G. Levitt \cite{levit} and has been investigated in depth by
D. Gaboriau \cite{gabor}. One of the major open problems here is the Fixed
Price problem, that is, whether every free action of a countable group has the
same cost. A. Kechris showed that for finitely generated groups and free
actions, the cost is monotonic with respect to weak containment
\cite[Corollary 10.14]{kechbook}, so Theorem \ref{fotetel} leads to the following.

\begin{corollary}
\label{maxcost}Let $\Gamma$ be a finitely generated group. Then, among free
p.m.p. actions of $\Gamma$, the cost is maximal on free factors of i.i.d.-s.
In particular, all free factors of i.i.d.-s have the same cost.
\end{corollary}

Since the cost of any free action of an infinite group is at least $1$, we get
that an infinite group $\Gamma$ has fixed price $1$ if and only if one of its
Bernoulli actions have cost $1$. We will generalize the monotonicity result of
A. Kechris to actions that are not necessarily free, but for that we need to
use an extension of the notion of cost called the \emph{groupoid cost}. For
details see Section \ref{cost}.

In the language of weak containment, one can express the well-known notion of
strong ergodicity as follows. An action of $\Gamma$ is \emph{strongly
ergodic}, if it is ergodic and it does not weakly contain the identity action
of $\Gamma$ on two points with equal mass.

In the presence of ergodicity, the absence of strong ergodicity forces various
structural restrictions on the action. An example of this is the theorem of V.
Jones and K. Schmidt \cite{jonessmith} that says that the equivalence relation
defined by such an action factors onto a nontrivial hyperfinite equivalence
relation. We will give another kind of restriction in the language of weak containment.

For a measure preserving action $f$ on $(X,\mathcal{B},\mu)$ let $f\times I$
denote the diagonal action of $\Gamma$ on $(X,\mathcal{B},\mu)\times
\lbrack0,1]$ where $\Gamma$ acts on the second coordinate trivially.

\begin{theorem}
\label{interval}Let $f$ be an ergodic p.m.p. action of a countable group. Then
the following are equivalent: \newline1) $f$ is not strongly ergodic;
\newline2) $f$ is weakly equivalent to $f\times I$.
\end{theorem}

This leads to the following convexity result, that can be viewed as a
`relative Glasner-Weiss theorem'.

For a compact topological space $K$ let $K^{\Gamma}$ denote the shift of
$\Gamma$ with base set $K$. This is a continuous action on a compact space.
Let $M(\Gamma,K)$ denote the set of $\Gamma$-invariant Borel measures on
$K^{\Gamma}$, endowed with the weak* topology and let $E(\Gamma,K)\subseteq
M(\Gamma,K)$ denote the set of ergodic measures. It is a standard fact that
$M(\Gamma,K)$ is a simplex where the set of extreme points equals
$E(\Gamma,K)$.

For an ergodic p.m.p. action $f$ on $(X,\mathcal{B},\mu)$ any Borel map
$\phi:X\rightarrow K$ defines a $\Gamma$-equivariant map $\Phi:X\rightarrow
K^{\Gamma}$ by setting, for $x\in X$ and $\gamma\in\Gamma$
\begin{equation}
\Phi(x)(\gamma)=\phi(f_{\gamma}(x))\text{.}\tag{1}%
\end{equation}
The measure $\Phi\circ\mu$ is an invariant measure on $K^{\Gamma}$.

Let
\[
E(f,K)=\left\{  \Phi\circ\mu\mid\phi\text{ a Borel map from }X\text{ to
}K\right\}  \text{.}%
\]
One can also describe $E(f,K)$ as the set of invariant measures on $K^{\Gamma
}$ that are factors of $f$. We have the following dichotomy, that can be
viewed as a `relative Glasner-Weiss theorem'.

\begin{theorem}
\label{convex}Let $f$ be an ergodic p.m.p. action of the countable group
$\Gamma$. If $f$ is strongly ergodic, then $\overline{E(f,K)}\subseteq
E(\Gamma,K)$. If $f$ is not strongly ergodic, then $\overline{E(f,K)}$ is convex.
\end{theorem}

The Glasner-Weiss theorem \cite{glasnweiss} says that when $\Gamma$ has
property (T), then $E(\Gamma,K)$ is closed, and when it does not have property
(T), then $\overline{E(\Gamma,K)}=M(\Gamma,K)$. We can deduce this result from
Theorem \ref{convex} as follows. By the Glasner-Thouvenot-Weiss theorem
\cite{glatuwei}, for every countable group $\Gamma$ there exists an ergodic
p.m.p. action $f$ of $\Gamma$ that weakly contains all ergodic actions of
$\Gamma$. This means that $E(\Gamma,K)\subseteq\overline{E(f,K)}$.

If $\Gamma$ has property (T), then by K. Schmidt \cite{schmidt2}, every
ergodic action of $\Gamma$ is strongly ergodic, in particular, $f$ is strongly
ergodic. So by Theorem \ref{convex} $\overline{E(f,K)}=E(\Gamma,K)$ and hence
$E(\Gamma,K)$ is closed.

If $\Gamma$ does not have property (T), then by the Connes-Weiss theorem,
$\Gamma$ has an ergodic action that is not strongly ergodic. But that action
is weakly contained in $f$, so $f$ is also not strongly ergodic. So by Theorem
\ref{convex}, $\overline{E(f,K)}$ is convex and since it contains all of
$E(\Gamma,K)$, it must equal $M(\Gamma,K)$, the convex hull of $E(\Gamma,K)$.
\bigskip

The paper is organized as follows. In Section \ref{Bernoulli} we prove Theorem
\ref{fotetel}. In Section \ref{cost} we define the groupoid cost and show that
for finitely generated groups, it is monotonic with respect to weak
containment. Finally, in Section \ref{convexity} we prove Theorem
\ref{interval} and Theorem \ref{convex}. \bigskip

\noindent\textbf{Acknowledgement.} The authors thank the Glasner family (Eli
and Yair) for establishing a link between them, that led to this paper.

\section{Weak containment of i.i.d. actions \label{Bernoulli}}

In this section we prove Theorem \ref{fotetel}. Throughout the paper, when it
is convenient, we will assume that the p.m.p. action is continuous on a
compact metric space. This can be easily seen by taking a Borel bijection
$\phi$ between $X$ and the Cantor set $K$ and then considering the shift
action on $K^{\Gamma}$ with the invariant measure $\Phi\circ\mu$.

For positive numbers $p_{1},\ldots,p_{d}$ with
\[%
{\displaystyle\sum\limits_{i=1}^{d}}
p_{i}=1
\]
let $\kappa(p_{1},\ldots,p_{d})$ denote the probability space on $\left\{
1,\ldots,d\right\}  $ where $P(i)=p_{i}$.

For any $T\subseteq\Gamma$ there is a natural projection $\pi_{T}$ from
$\left\{  1,\ldots,d\right\}  ^{\Gamma}$ to $\left\{  1,\ldots,d\right\}
^{T}$. The sets of the form $\pi_{T}^{-1}(D)$ with $D\subseteq\left\{
1,\ldots,d\right\}  ^{T}$ are called the \emph{cylinder sets with respect to
}$T$.

\begin{lemma}
\label{reduction}Let $f$ be a p.m.p. action of the countable group $\Gamma$ on
$(X,\mathcal{B},\mu)$. Then $f$ weakly contains $\kappa(p_{1},\ldots
,p_{d})^{\Gamma}$ if for all finite subsets $F$ of $\Gamma$ and $\delta>0$
there exists a Borel map $\phi:X\rightarrow\left\{  1,\ldots,d\right\}  $ such
that for all $\alpha\in\left\{  1,\ldots,d\right\}  ^{F}$
\[
\left\vert \mu\left(  \left\{  x\in X\mid\phi(f_{\gamma}x)=\alpha
(\gamma)\text{ for all }\gamma\in F\right\}  \right)  -%
{\displaystyle\prod\limits_{\gamma\in F}}
p_{\alpha(\gamma)}\right\vert <\delta\text{.}%
\]

\end{lemma}

\noindent\textbf{Proof.} Assume that the condition holds. Let $g$ denote the
shift action of $\Gamma$ on $\kappa(p_{1},\ldots,p_{d})^{\Gamma}$ and let
$\lambda$ denote the measure on $\kappa(p_{1},\ldots,p_{d})^{\Gamma}$.

Let $Y_{1},\ldots,Y_{n}\subseteq\left\{  1,\ldots,d\right\}  ^{\Gamma}$ be
Borel subsets, let $S\subseteq\Gamma$ be a finite symmetric set containing the
identity and let $\varepsilon>0$.

Then for any $r>0$, there exists a finite subset $T$ of $\Gamma$, and subsets
$D_{i}\subseteq\left\{  1,\ldots,d\right\}  ^{T}$ such that for all $1\leq
i\leq n$ we have
\[
\lambda(Y_{i}\vartriangle\pi_{T}^{-1}(D_{i}))<r\text{.}%
\]
Here $\vartriangle$ denotes symmetric difference. We now apply the condition
of the lemma for $F=TS$ and a $\delta>0$ to be chosen later.

For $1\leq i\leq n$ let
\[
X_{i}=\Phi^{-1}(\pi_{T}^{-1}(D_{i}))
\]
where $\Phi$ is the map defined by 1) in the Introduction. Then for
sufficiently small $r$ and $\delta$ it follows that%

\[
\left\vert \mu(f_{\gamma}X_{i}\cap X_{j})-\lambda(g_{\gamma}Y_{i}\cap
Y_{j})\right\vert <\varepsilon\text{ \ (}1\leq i,j\leq n,\gamma\in S\text{).}%
\]

The lemma is proved. $\square$

\begin{lemma}
\label{free}Let $f$ be a free p.m.p. action of the countable group $\Gamma$ on
the compact metric space $(X,d,\mu)$. Then for any finite symmetric subset $F$
of $\Gamma$ and any $\varepsilon>0$ there exists $s>0$ such that
\begin{equation}
\mu\left(  \left\{  x\in X\mid d(f_{\gamma}x,f_{\gamma^{\prime}}x)>s\text{ for
all }\gamma\neq\gamma^{\prime}\in F\right\}  \right)  >1-\varepsilon\tag{*}%
\end{equation}
and
\begin{equation}
\mu\times\mu\left(  \left\{  (x,x^{\prime})\in X\times X\mid d(f_{\gamma
}x,f_{\gamma^{\prime}}x^{\prime})>s\text{ for all }\gamma,\gamma^{\prime}\in
F\right\}  \right)  >1-\varepsilon\tag{**}%
\end{equation}

\end{lemma}

In other terms, for this $s$, most of the $F$-neighbourhoods of pairs of
points in $X$ are $s$-apart. We leave the details to the reader. \bigskip

We are ready to prove Theorem \ref{fotetel}. \bigskip

\noindent\textbf{Proof of Theorem \ref{fotetel}.} We will verify that the
conditions of Lemma \ref{reduction} hold. Let $F$ be a finite subset of
$\Gamma$ and let $\delta>0$. Let $\varepsilon>0$ be determined later. Use
Lemma \ref{free} with this $\varepsilon$ to find $s$ that satisfies the
conclusions of the lemma. Let $\left\{  B_{j}\mid1\leq j\leq J\right\}  $ be a
Borel partition of $X$ into sets of diameter less than $s$. Let $U_{1}%
,\ldots,U_{J}$ be independent random variables with distribution $\kappa
(p_{1},\ldots,p_{d})$. Define a random map $\phi:X\rightarrow\left\{
1,\ldots,d\right\}  $ by setting
\[
\phi(x)=U_{j}\text{ when }x\in B_{j}\text{.}%
\]
We shall show that with high probability, this map will satisfy the conditions
of Lemma \ref{reduction}. In particular, this will imply the existence of such
a map.

We shall use the so-called second moment method. For $\alpha\in\left\{
1,\ldots,d\right\}  ^{F}$ let
\[
G_{\alpha}=\left\{  x\in X\mid\phi(f_{\gamma}x)=\alpha(\gamma)\text{ for all
}\gamma\in F\right\}  \text{.}%
\]
Now $\mu(G_{\alpha})$ is a random variable with expected value given by
\[
\boldsymbol{E}\mu(G_{\alpha})=\boldsymbol{E}%
{\displaystyle\int\limits_{X}}
1_{G_{\alpha}}(x)d\mu(x)=%
{\displaystyle\int\limits_{X}}
\boldsymbol{E}1_{G_{\alpha}}(x)d\mu(x)\text{.}%
\]
For $x\in X$ such that the $F$-neighbourhood of $x$ is $s$-apart (that is, the
set defined in (*) of Lemma \ref{free}) the expected value $\boldsymbol{E}%
1_{G_{\alpha}}(x)=%
{\displaystyle\prod\limits_{\gamma\in F}}
p_{\alpha(\gamma)}$ because the $\phi$-values at the neighbours are
independent. It follows that
\[
\left\vert \boldsymbol{E}\mu(G_{\alpha})-%
{\displaystyle\prod\limits_{\gamma\in F}}
p_{\alpha(\gamma)}\right\vert <2\varepsilon\text{.}%
\]
We shall now estimate the variance of $\mu(G_{\alpha})$. We have
\begin{align*}
\boldsymbol{E}\mu^{2}(G_{\alpha}) &  =\boldsymbol{E}%
{\displaystyle\int\limits_{X}}
1_{G_{\alpha}}(x)d\mu(x)%
{\displaystyle\int\limits_{X}}
1_{G_{\alpha}}(x^{\prime})d\mu(x^{\prime})=\\
&  =\boldsymbol{E}%
{\displaystyle\int\limits_{X}}
{\displaystyle\int\limits_{X}}
1_{G_{\alpha}}(x)1_{G_{\alpha}}(x^{\prime})d(\mu\times\mu)(x,x^{\prime})=\\
&  =%
{\displaystyle\int\limits_{X}}
{\displaystyle\int\limits_{X}}
\boldsymbol{E}\mathbb{(}1_{G_{\alpha}}(x)1_{G_{\alpha}}(x^{\prime}%
))d(\mu\times\mu)(x,x^{\prime})
\end{align*}
For $x,x^{\prime}\in X$ such that the $F$-neighbourhoods of $x$ and
$x^{\prime}$ are $s$-apart (that is, when both $x$, $x^{\prime}$ satisfy (*)
and $(x,x^{\prime})$ satisfies (**) of Lemma \ref{free}) then%
\[
\boldsymbol{E}1_{G_{\alpha}}(x)1_{G_{\alpha}}(x^{\prime})=\left(
{\displaystyle\prod\limits_{\gamma\in F}}
p_{\alpha(\gamma)}\right)  ^{2}%
\]
because the $\phi$-values at the neighbours are independent. The $\mu\times
\mu$-measure of these pairs of points is at least $1-3\varepsilon$. It follows
that
\[
\mathrm{Var}(\mu(G_{\alpha}))=\boldsymbol{E}\mu^{2}(G_{\alpha}%
)-(\boldsymbol{E}\mu(G_{\alpha}))^{2}<12\varepsilon\text{.}%
\]
Now Chebyshev's inequality says that the probability
\[
P(\left\vert \mu(G_{\alpha})-\boldsymbol{E}\mu(G_{\alpha})\right\vert
>c)\leq\frac{\mathrm{Var}(\mu(G_{\alpha}))}{c^{2}}%
\]
Setting $c=\varepsilon^{1/3}$ we get
\[
P(\left\vert \mu(G_{\alpha})-\boldsymbol{E}\mu(G_{\alpha})\right\vert
>\varepsilon^{1/3})\leq12\varepsilon^{1/3}%
\]
Setting $\varepsilon$ small enough, we get that $\left\vert \mu(G_{\alpha
})-\boldsymbol{E}\mu(G_{\alpha})\right\vert $ is arbitrarily small for all
$\alpha\in\left\{  1,\ldots,d\right\}  ^{F}$ simultaneously. Hence the
conditions of Lemma \ref{reduction} hold and $f$ weakly contains $\kappa
(p_{1},\ldots,p_{d})^{\Gamma}$ for any finite distribution $\kappa
(p_{1},\ldots,p_{d})$.

Since any probability distribution can be approximated by finite
distributions, we get that $f$ weakly contains all Bernoulli actions and the
theorem has been established. $\square$

\section{Groupoid cost and weak containment \label{cost}}

In this section we will introduce the groupoid cost of a p.m.p. action and
show that for finitely generated groups, it is monotonic with respect to weak
containment. Groupoid cost has been introduced in \cite{abnikheg} as an
extension of cost for non-free p.m.p. actions. The monotonicity result was
proved by A. Kechris \cite[Corollary 10.14]{kechbook} for the original cost
and free actions.

Let $f$ be a p.m.p. action of the countable group $\Gamma$ on $(X,\mathcal{B}%
,\mu)$. Without loss of generality, we can assume that the space $X$ is
compact. Endow $\Gamma$ with the discrete topology and the counting measure.
Consider $X\times\Gamma$, endowed with the product Borel structure and the
product measure $\widetilde{\mu}$. We define a partial product on
$X\times\Gamma$ as follows. The product of $(x_{1},\gamma_{1})$ and
$(x_{2},\gamma_{2})$ is only defined, when $f_{\gamma_{1}}x_{1}=x_{2}$; in
this case let $(x_{1},\gamma_{1})\cdot(x_{2},\gamma_{2})=(x_{1},\gamma
_{1}\gamma_{2})$. The inverse of $(x,\gamma)$ is defined as $(f_{\gamma
}x,\gamma^{-1})$ and is denoted by $(x,\gamma)^{-1}$, so we have
$(x,\gamma)\cdot(x,\gamma)^{-1}=(x,e)$. Then $X\times\Gamma$ forms a groupoid
with respect to this partial product.

Now we define the \emph{subset groupoid} $\mathcal{G}_{f}$ as follows.
Elements of $\mathcal{G}_{f}$ are Borel subsets of $X\times\Gamma$. For
$A,B\in\mathcal{G}_{f}$ let
\[
A\cdot B=\left\{  a\cdot b\mid a\in A\text{, }b\in B\text{ when }a\cdot
b\text{ is defined}\right\}  \text{.}%
\]

Let $E=X\times\{e\}$ where $e$ is the identity element of $\Gamma$. Then for
all $A\in\mathcal{G}_{f}$ we have $AE=EA=A$. When $f$ is a continuous action
on a compact space, the partial product is continuous as well, which in
particular implies that the product of open elements of $\mathcal{G}_{f}$ is
open and the product of compact elements is compact.

We say that $A\in\mathcal{G}_{f}$ \emph{generates the action }$f$, if
\[%
{\displaystyle\bigcup\limits_{n=1}^{\infty}}
\left(  A\cup A^{-1}\cup E\right)  ^{n}=X\times\Gamma\text{.}%
\]
In particular, for any generating set $S$ of $\Gamma$, $X_{S}=X\times S$
generates $f$. Note that this is a topological condition and is independent of
the measure.

The \emph{groupoid cost} of $f$ is defined as
\[
\mathrm{gcost}(f)=\inf_{A\text{ generates }f}\widetilde{\mu}(A)\text{.}%
\]
When $f$ is a free action, $X\times\Gamma$ can be identified with the
equivalence relation generated by $f$ on $X$, elements of $\mathcal{G}_{f}$
are graphings in this relation and $A\in\mathcal{G}_{f}$ generates $f$ if and
only if the correspoding graphing generates the equivalence relation. Hence,
for free actions, $\mathrm{gcost}(f)=\mathrm{cost}(f)$ and so the groupoid
cost is indeed an extension of cost to non-free actions. Note that the
original notion of $\mathrm{cost}$ required that the graphing generates the
relation only up to a nullset, but it is easy to see that without loss of
generality, one can assume full generation.

Trivially, for every generator $A\in\mathcal{G}_{f}$ and $\varepsilon>0$ there
exists an open generator $B\in\mathcal{G}_{f}$ such that $A\subseteq B $ and
$\widetilde{\mu}(B)\leq\widetilde{\mu}(A)+\varepsilon$. So, in the definition
of gcost, it is enough to consider open generators. The following lemma is
essentially contained in \cite{abnikheg}.

\begin{lemma}
\label{clopen}Assume that $X$ is compact, totally disconnected and $\Gamma$ is
finitely generated. Then
\[
\mathrm{gcost}(f)=\inf\widetilde{\mu}(A)
\]
where $A$ varies over compact, clopen generators for $f$.
\end{lemma}

\noindent\textbf{Proof.} Let $\Gamma$ be generated by the finite set $S$. Let
$\mathcal{C}$ be a countable clopen basis for the topology of $X$. List the
elements of $\mathcal{C}$ as $C_{1},C_{2},\ldots$ For an open set
$O\in\mathcal{B}$ and $n>0$ let
\[
O_{n}=%
{\displaystyle\bigcup\limits_{\substack{1\leq i\leq n \\C_{i}\subseteq O}}}
C_{i}.
\]
Let $A\in\mathcal{G}_{f}$ be an open generating element for $\mathcal{G}_{f} $
with $\widetilde{\mu}(A)<\infty$. Then $A=\cup_{\gamma\in\Gamma}%
A(\gamma)\times\{\gamma\}$ where $A(\gamma)\times\{\gamma\}=A\cap\left(
X\times\{\gamma\}\right)  $ and $A(\gamma)$ is open. For $n>0$ let $A_{n}%
=\cup_{\gamma\in\Gamma}A(\gamma)_{n}\times\{\gamma\}$. Then $A_{n}$ is compact
(since $\widetilde{\mu}(A)<\infty$) and clopen. Let $B_{n}=A_{n}\cup
A_{n}^{-1}\cup E$. Then
\[%
{\displaystyle\bigcup\limits_{n=1}^{\infty}}
B_{n}^{n}=X\times\Gamma
\]
since every possible finite product in $A\cup A^{-1}\cup E$ is realized in
$B_{n}$ for large enough $n$. In particular,
\[%
{\displaystyle\bigcup\limits_{n=1}^{\infty}}
\left(  B_{n}^{n}\cap X_{S}\right)  =X_{S}%
\]
but since $B_{n}^{n}\cap X_{S}$ is open and $X_{S}$ is compact, there exists
an $n>0$ with $X_{S}\subseteq B_{n}^{n}$. Since $X_{S}$ generates $f$, we get
that already $A_{n}$ generates $f$.

It follows that any open generator for $f$ with finite measure contains a
compact, clopen generator. The lemma is proved. $\square$\bigskip

For the convenience of the reader, we now recall some definitions from the
Introduction. For a compact topological space $K$ let $K^{\Gamma}$ denote the
shift of $\Gamma$ with base set $K$. This is a continuous action on a compact
space. Let $M(\Gamma,K)$ denote the set of $\Gamma$-invariant Borel measures
on $K^{\Gamma}$, endowed with the weak* topology and let $E(\Gamma,K)\subseteq
M(\Gamma,K)$ denote the set of ergodic measures.

For an ergodic p.m.p. action $f$ on $(X,\mathcal{B},\mu)$ any Borel map
$\phi:X\rightarrow K$ defines a $\Gamma$-equivariant map $\Phi:X\rightarrow
K^{\Gamma}$ by setting, for $x\in X$ and $\gamma\in\Gamma$
\[
\Phi(x)(\gamma)=\phi(f_{\gamma}(x))\text{.}%
\]
The measure $\Phi\circ\mu$ is an invariant measure on $K^{\Gamma}$.

Let
\[
E(f,K)=\left\{  \Phi\circ\mu\mid\phi\text{ a Borel map from }X\text{ to
}K\right\}  \text{.}%
\]
One can also describe $E(f,K)$ as the set of invariant measures on $K^{\Gamma
}$ that are factors of $f$.

\begin{lemma}
\label{closure}Let $f$ and $g$ be p.m.p. actions of the countable group
$\Gamma$ on $(X,\mathcal{B},\mu)$ and $(Y,\mathcal{C},\nu)$, respectively.
Then $f$ weakly contains $g$ if and only if
\begin{equation}
\overline{E(f,K)}\supseteq E(g,K)\tag{*}%
\end{equation}
for any compact topological space $K$. The same result holds if we only
consider the class of finite spaces.
\end{lemma}

\noindent\textbf{Proof.} The (*) condition for all finite spaces $K$ is easily
seen to be equivalent with the original definition of weak containment.

If $L$ is a finite subset of $K$, then $L^{\Gamma}$ is contained in
$K^{\Gamma}$. The invariant measures concentrated on $L^{\Gamma}$ as $L$
varies over the finite subsets of $K$, are dense in $M(\Gamma,K)$.

Moreover, it follows that
\[
\overline{%
{\displaystyle\bigcup\limits_{L\subseteq K\text{ finite}}}
E(g,L)}\supseteq E(g,K)\text{.}%
\]
This is because
\[
E(g,K)=\left\{  \Phi\circ\nu\mid\phi\text{ a Borel map from }Y\text{ to
}K\right\}  \text{.}%
\]
and the Borel map $\phi$ can be well approximated by maps $\phi_{n}$ with
finite range $L_{n}$ in such a way that $\Phi_{n}\circ\nu$ converges to
$\Phi\circ\nu$. This yields
\[
E(g,K)\subseteq\overline{%
{\displaystyle\bigcup\limits_{L\subseteq K\text{ finite}}}
E(g,L)}\subseteq\overline{%
{\displaystyle\bigcup\limits_{L\subseteq K\text{ finite}}}
E(f,L)}\subseteq\overline{E(f,K)}\text{.}%
\]

The lemma is proved. $\square$

\begin{theorem}
\label{monoton}Let $f$ and $g$ be p.m.p. actions of the finitely generated
group $\Gamma$ on $(X,\mathcal{B},\mu)$ and $(Y,\mathcal{C},\nu)$,
respectively, such that $f$ weakly contains $g$. Then $\mathrm{gcost}%
(f)\leq\mathrm{gcost}(g)$.
\end{theorem}

\noindent\textbf{Proof.} Let $K$ be the standard Cantor set. Let $\phi$ be a
Borel isomorphism from $Y$ to $K$. Then the shift action on $K^{\Gamma}$ with
the invariant measure $\nu^{\prime}=\Phi\circ\nu$ is isomorphic to $g$. Let
$h$ be the shift action of $\Gamma$ on $K^{\Gamma}$ with the invariant measure
$\nu^{\prime}$. Clearly, $\mathrm{gcost}(h)=\mathrm{gcost}(g)$. Using Lemma
\ref{clopen}, for any $\varepsilon>0$ there exists a clopen, compact
$A\in\mathcal{G}_{h}$ that generates $h$ and $\widetilde{\nu}^{\prime
}(A)<\mathrm{gcost}(h)+\varepsilon$. Since weak convergence of probability
measures $\rho_{n}$ on $K^{\Gamma}$ implies convergence of $\widetilde{\rho
_{n}}$ on clopen compact subsets of $K^{\Gamma}\times\Gamma$, using Lemma
\ref{closure} there exists a map $\psi:X\rightarrow K$ such that
\[
\widetilde{\Psi\circ\mu}(A)\leq\widetilde{\nu}^{\prime}(A)+\varepsilon\text{.}%
\]
The condition of generating is independent of the invariant measure, so $A$
generates the shift action on $K^{\Gamma}$ with the invariant measure
$\Psi\circ\mu$ as well. This action is a factor of $f$, and hence its gcost is
an upper bound for the gcost of $f$. This implies
\[
\mathrm{gcost}(f)\leq\widetilde{\Psi\circ\mu}(A)\leq\widetilde{\nu}^{\prime
}(A)+\varepsilon<\mathrm{gcost}(g)+2\varepsilon.
\]

The theorem holds. $\square$\bigskip

\noindent\textbf{Proof of Corollary \ref{maxcost}.} Let $\Gamma$ be a finitely
generated group and let $f$ be a free p.m.p. action of $\Gamma$. Then, using
Theorem \ref{fotetel}, $f$ weakly contains all Bernoulli actions, so by
Theorem \ref{monoton} the gcost of $f$ is at most the gcost of any Bernoulli
action of $\Gamma$, in particular, the cost of $f$ is at most the cost of any
free Bernoulli action of $\Gamma$. $\square$\bigskip

Note that any nontrivial Bernoulli action of an infinite group is free. Also,
it is easy to see that when $\Gamma$ is torsion free, then every nontrivial
factor of i.i.d. is free. However, in the presence of torsion, this last
statement fails to hold.

\section{Ergodic actions and convexity \label{convexity}}

Let $f$ be a p.m.p. action of the countable group $\Gamma$ on $(X,\mathcal{B}%
,\mu)$. A sequence of subsets $A_{n}\in\mathcal{B}$ is \emph{almost
invariant}, if
\[
\lim_{n\rightarrow\infty}\mu(f_{\gamma}A_{n}\vartriangle A_{n})=0\text{ for
all }\gamma\in\Gamma\text{.}%
\]
The sequence is \emph{trivial}, if $\lim_{n\rightarrow\infty}\mu(A_{n}%
)(1-\mu(A_{n}))=0$. We say that the action $f$ is \emph{strongly ergodic}, if
every almost invariant sequence is trivial.

The following lemma is due to K. Schmidt \cite{schmidt}.

\begin{lemma}
[Schmidt]Let $f$ be a p.m.p. action of the countable group $\Gamma$ on
$(X,\mathcal{B},\mu)$, such that $f$ is ergodic, but not strongly ergodic.
Then for all $0<\lambda<1$ there exists an almost invariant sequence $A_{n}$
such that $\mu(A_{n})=\lambda$ ($n\geq0$).
\end{lemma}

In particular, the action is strongly ergodic if and only if it is ergodic and
does not weakly contain the identity action of $\Gamma$ on two points with
equal mass.

The following lemma forms a part of the proof of Schmidt's lemma. For
completeness, we provide a short proof.

\begin{lemma}
\label{schmidt}Let $f$ be an ergodic p.m.p. action of the countable group
$\Gamma$ on the compact space $(X,\mathcal{B},\mu)$, let $0<\lambda<1$ and let
$A_{n}$ be an almost invariant sequence such that
\[
\lim_{n\rightarrow\infty}\mu(A_{n})=\lambda\text{.}%
\]
Let $\mu_{n}$ be the measure on $X$ defined by
\[
\mu_{n}(Y)=\mu(A_{n}\cap Y)\text{ \ (}Y\in\mathcal{B}\text{).}%
\]
Then $\mu_{n}$ weak* converges to $\lambda\mu$.
\end{lemma}

\noindent\textbf{Proof.} Let $n_{k}$ be a sequence such that $\mu_{n_{k}}$ is
weakly convergent and let $\nu$ be the limit measure. Then $\nu$ is invariant
under $\Gamma$ and is bounded by $\mu/\lambda$, so by ergodicity, $\nu$ is a
scalar multiple of $\mu$. But $\nu(X)=\lim\mu(A_{n})=\lambda$ which proves the
lemma. $\square$\bigskip

\noindent\textbf{Proof of Theorem \ref{interval}.} Let $f$ be an ergodic
p.m.p. action of the countable group $\Gamma$ on $(X,\mathcal{B},\mu)$.

Assume that $f$ is weakly equivalent to $f\times I$. Since $f\times I$ factors
on the identity action $b$ of $\Gamma$ on two points with equal mass, this
means that $f$ weakly contains $b$. Hence $f$ can not be strongly ergodic.
\newline

Now assume $f$ is not strongly ergodic. We claim that $f$ weakly contains
$f\times b$. We can assume that the space $(X,\mathcal{B},\mu)$ is totally
disconnected and compact and $\Gamma$ acts by continuous maps. Let
$\mathcal{C}$ denote the set of clopen sets in $\mathcal{B}$. Then
$\mathcal{C}$ is a finitely additive algebra and for any $B\in\mathcal{B}$ and
$\varepsilon>0$ there exists $C\in\mathcal{C}$ with $\mu(B\vartriangle
C)<\varepsilon$.

Hence, it is enough to check the weak containment condition for clopen
subsets. In addition, by Schmidt's lemma, there exists an almost invariant
sequence $A_{n}\in\mathcal{C}$ with
\[
\lim_{n\rightarrow\infty}\mu(A_{n})=1/2\text{.}%
\]
Lemma \ref{schmidt} implies that for every clopen set $C\in\mathcal{C}$, we
have
\[
\lim_{n\rightarrow\infty}\mu(C\cap A_{n})=\mu(C)/2\text{.}%
\]
Indeed, weak* convergence of measures implies convergence on clopen sets.

Let $Y_{1}^{\prime},\ldots,Y_{k}^{\prime},Y_{1}^{\prime\prime},\ldots
,Y_{k}^{\prime\prime}\in\mathcal{C}$, let $S\subseteq\Gamma$ be a finite
subset and let $\varepsilon>0$. Then by the almost invariance of $A_{n}$, for
all $\gamma\in S$ we have
\[
\lim_{n\rightarrow\infty}\mu(f_{\gamma}A_{n}\vartriangle A_{n})=0
\]

Also for all $1\leq i,j\leq k$ and $\gamma\in S$ we have
\[
\lim_{n\rightarrow\infty}\mu(A_{n}\cap f_{\gamma}Y_{i}^{\prime}\cap
Y_{j}^{\prime})-\mu(f_{\gamma}Y_{i}^{\prime}\cap Y_{j}^{\prime})/2=0
\]
and
\[
\lim_{n\rightarrow\infty}\mu(A_{n}\cap f_{\gamma}Y_{i}^{\prime\prime}\cap
Y_{j}^{\prime\prime})-\mu(f_{\gamma}Y_{i}^{\prime\prime}\cap Y_{j}%
^{\prime\prime})/2=0
\]
Let $n$ be large enough such that all the above quantities (listed in the
above three displays) are less than $\varepsilon$. Now for $1\leq i\leq k$
let
\[
X_{i}=\left(  A_{n}\cap Y_{i}^{\prime}\right)  \cup\left(  A_{n}^{c}\cap
Y_{i}^{\prime\prime}\right)
\]
where $A_{n}^{c}$ denotes the complement of $A_{n}$.

Now let $Y_{1},\ldots,Y_{k}$ be Borel subsets of $X\times\{0,1\}$. Let
\[
Y_{i}^{\prime}\times\{0\}=Y_{i}\cap(X\times\{0\})
\]
and let
\[
Y_{i}^{\prime\prime}\times\{1\}=Y_{i}\cap(X\times\{1\})\text{.}%
\]
Applying the above argument for $Y_{i}^{\prime}$ and $Y_{i}^{\prime\prime} $
we obtain the desired measurable subsets $X_{1},\ldots,X_{n}$. The claim holds
and $f$ weakly contains $f\times b$.

Iterating the argument, we get that $f$ weakly contains $f\times b_{2^{n}}$
where $b_{2^{n}}$ is the trivial action of $\Gamma$ on $\left\{
1,\ldots,2^{n}\right\}  $ with uniform measure. An easy approximation argument
yields that $f$ weakly contains $f\times I$. Since $f$ is a factor of $f\times
I$, we get that $f$ is weakly equivalent to $f\times I$.

The theorem is proved. $\square$\bigskip

\noindent\textbf{Proof of Theorem \ref{convex}.} Let $f$ be an ergodic p.m.p.
action of the countable group $\Gamma$.

Assume that $f$ is strongly ergodic. Assume by contradiction that there exists
a sequence of Borel maps $\phi_{n}:X\rightarrow K$ such that $\Phi_{n}\circ
\mu$ weak* converges to a non-ergodic measure $\lambda$ on $K^{\Gamma}$. Let
$A$ be a $\Gamma$-invariant Borel set with $0<\lambda(A)<1$. Let us
approximate $A$ with a clopen sequence $A_{n}$ with respect to $\lambda$. It
is easy to see that $A_{n}$ is almost invariant with respect to $\lambda$.
Then there exists $k(n)$ such that $\Phi_{k(n)}^{-1}(A_{n})$ gives us a
nontrivial almost invariant sequence for $f$, a contradiction. Hence
$\overline{E(f,K)}\subseteq E(\Gamma,K)$.

Assume that $f$ is not strongly ergodic. Let $\mu^{\prime}$ be the product
measure $\mu\times\{1/2,1/2\}$. Let $\phi_{0},\phi_{1}:X\rightarrow K$ be
Borel maps and let $\phi:X\times\{0,1\}\rightarrow K$ be defined by
\[
\phi((x,i))=\phi_{i}(x)\text{ for }i=0,1\text{.}%
\]
Then by Theorem \ref{interval}, $f$ weakly contains $f\times b$, so by Lemma
\ref{convex}, $\Phi\circ\mu^{\prime}$ is in $\overline{E(f,K)}$. But
\[
\Phi\circ\mu^{\prime}=\frac{1}{2}\Phi_{0}\circ\mu+\frac{1}{2}\Phi_{1}\circ
\mu\text{.}%
\]
This gives that for $\lambda,\lambda^{\prime}\in E(f,K)$, $\frac{1}{2}%
\lambda+\frac{1}{2}\lambda^{\prime}\in\overline{E(f,K)}$. Hence $\overline
{E(f,K)}$ is convex.

The theorem holds. $\square$

\end{document}